\documentclass[titlepage,11pt]{article}
\usepackage{latexsym}
\usepackage{amsfonts}
\usepackage{amsmath}
\newcommand\blackslug{\hbox{\hskip 1pt \vrule width 4pt height 8pt depth 1.5pt
        \hskip 1pt}}
\newcommand\bbox{\hfill \quad \blackslug \bigbreak}
\def\d{\hbox{-}}

\title{Towards Erd\H{o}s-Hajnal for graphs with no 5-hole}
\author{Maria Chudnovsky\thanks{Supported by NSF grant DMS-1550991.
This material is based upon work supported in part by the U. S. Army
Research Laboratory and the U. S. Army Research Office under grant
number
W911NF1610404.}\\
Princeton University, Princeton, NJ 08544
\\
\\
Jacob Fox\thanks{Supported by a Packard Fellowship and NSF Career Award DMS-1352121.}\\
Stanford University, Stanford, CA 94305-2125
\\
\\
Alex Scott\thanks{Supported by a Leverhulme Research
Fellowship.}\\
Mathematical Institute, University of Oxford, Oxford OX2 6GG, UK
\\
\\
Paul Seymour\thanks{Supported by ONR grant N00014-14-1-0084 and NSF
grant DMS-1265563.}\\
Princeton University, Princeton, NJ 08544
\\
\\
Sophie Spirkl\\
Princeton University, Princeton, NJ 08544}

\date{March 3, 2018; revised \today}

\newtheorem{thm}{}[section]

\newcommand{\Proof}{\noindent{\bf Proof.}\ \ }

\begin{document}
\maketitle
\begin{abstract} 
The Erd\H{o}s-Hajnal conjecture says that for every graph $H$ there exists $c>0$
such that 
$$\max(\alpha(G),\omega(G))\ge n^c$$ 
for every $H$-free graph $G$ with $n$ vertices, and this is still open when $H=C_5$.
Until now the best bound known on $\max(\alpha(G),\omega(G))$ for $C_5$-free graphs 
was the general bound of Erd\H{o}s and Hajnal, that for all $H$,
$$\max(\alpha(G),\omega(G))\ge 2^{\Omega(\sqrt{\log n })}$$ if $G$ is $H$-free. We improve this
when $H=C_5$ to 
$$\max(\alpha(G),\omega(G))\ge 2^{\Omega(\sqrt{\log n \log \log n  })}.$$

\end{abstract}

\section{Introduction}

All graphs in this paper are finite and have no loops or parallel edges, and the cardinalities of the largest stable sets and 
cliques in a graph $G$
are denoted by $\alpha(G),\omega(G)$ respectively.
If $G,H$ are graphs,
we say that $G$ {\em contains} $H$ if some induced subgraph of $G$ is isomorphic to $H$, and $G$ is
{\em $H$-free} otherwise.

The Erd\H{o}s-Hajnal conjecture~\cite{EH0,EH} asserts:

\begin{thm}\label{EHconj}
{\bf Conjecture: }For every graph $H$, there exists $\epsilon>0$ such that every $H$-free graph $G$ satisfies
$$\max(\alpha(G),\omega(G))\ge |V(G)|^\epsilon.$$
\end{thm}
This is true for all $H$ with at most four vertices, but
is open when $H=C_5$  ($C_5$ denotes the cycle of length five).
The problem for $C_5$ has attracted a good deal of unsuccessful attention, for several reasons; not only is $C_5$ arguably
the smallest open case of \ref{EHconj}, but also it has a good amount of symmetry, and more importantly, 
by excluding $C_5$ we exclude its complement as well. (Excluding both a graph and its complement is an approach that
has been quite fruitful lately, for instance~\cite{bonamy, lagoutte}.) So we are happy to report some progress at last.

The best general bound for the Erd\H{o}s-Hajnal conjecture to date was proved by Erd\H{o}s and Hajnal in~\cite{EH},
namely:

\begin{thm}\label{EHpartial}
For every graph $H$, there exists $c>0$ such that 
$$\max(\alpha(G),\omega(G))\ge 2^{c\sqrt{\log n}}$$
for every $H$-free graph $G$ with $n>0$ vertices.
\end{thm}
(Logarithms are to base two, throughout the paper.)
Until now, this was also the best bound known when $H=C_5$, but in this paper we will improve it to:
\begin{thm}\label{mainthm}
There exists $c>0$ such that 
$$\max(\alpha(G),\omega(G))\ge 2^{c\sqrt{\log  n \log \log  n  }}$$
for every $C_5$-free graph $G$ with $n>1$ vertices.
\end{thm}

If $A,B\subseteq V(G)$ are disjoint and nonempty, the {\em edge-density} between them means the number of 
edges joining $A,B$,
divided by $|A|\cdot |B|$. The proof of \ref{mainthm} is via the following conjecture of
Conlon, Fox and Sudakov~\cite{fox}: 

\begin{thm}\label{foxdensity}
{\bf Conjecture: }For every graph $H$ there exist $\epsilon,\sigma > 0$ such that
for every $H$-free graph $G$ on $n>1$ vertices, and all $c$ with $0\le c\le 1/2$, $V(G)$ contains two
disjoint subsets $A,B$ with $|A|\ge \epsilon c^{\sigma} n$
and $|B|\ge \epsilon n$, such that the edge-density between $A,B$ is either at most $c$
or at least $1-c$.
\end{thm}

This has not been proved so far for any graph $H$ with more than four vertices, but in this paper we prove it
for $H=C_5$ (with $\sigma=1$), and this is the key to proving \ref{mainthm}.
We first prove it for sparse graphs $G$,  and then use a theorem of R\"odl to deduce it in 
general (both in the next section). The proof of \ref{mainthm} is completed in section 3.

We remark that \ref{foxdensity} (for all $H$) is equivalent to the same statement for sparse graphs (for all $H$),
because of the theorem of R\"odl discussed in the next section; but for sparse graphs we can prove \ref{foxdensity}
for many more graphs $H$ than just $C_5$ (for instance, for all bipartite $H$, and all cycles of length at least four). 
These results will appear in a later paper~\cite{polysparse}.
But $C_5$ is still the largest graph $H$
for which we can show that both $H$ and its complement satisfy \ref{foxdensity} in sparse graphs, and so the largest
for which we can prove \ref{foxdensity}. 

\section{Sparse graphs}

In this section we prove \ref{foxdensity} for $H=C_5$, and first we prove it when $G$ is sufficiently sparse. 
Let us say the {\em closed
degree} of a vertex is one more than its degree. (Counting cardinalities of subsets works out more conveniently using closed degree.)
For disjoint $A,B\subseteq V(G)$, we say $A$ is 
{\em anticomplete} to $B$ if there are no edges between $A$ and $B$.
We will prove:

\begin{thm}\label{get5hole}
For all $c$ with $0<c\le 1/2$,
and every graph $G$ with $n>0$ vertices, if $G$ satisfies:
\begin{itemize}
\item every vertex has closed degree at most $n/16$, and
\item for every two disjoint subsets $A,B\subseteq V(G)$ with $|A|\ge  c n/2$ and $|B|\ge n/16$,
the edge-density between $A,B$ is at least $c$,
\end{itemize}
then $G$ contains $C_5$.
\end{thm}
\Proof
Let $0<c\le 1/2$,  and let $G,n$ be as in the theorem.
Since every vertex has closed degree at most $n/16$, it follows that $n\ge 16$ and in particular,
$\lfloor n/2 \rfloor\ge n/4$.
Choose a set $N_0\subseteq V(G)$ of cardinality $\lfloor n/2 \rfloor$. It follows that $|N_0|\ge n/4\ge cn/2$, and
so the edge-density between $N_0$ and its complement is at least $c$.
In particular, some vertex in $N_0$ has at least $cn/2$ neighbours. 

Let $v_1$ be a vertex of degree at least $cn/2$, let $N_1$ be the set of all neighbours of $v_1$, and 
let $Z_2=V(G)\setminus (N_1\cup \{v_1\})$.
Since $|N_1|+1\le n/16$, it follows that $|Z_2|\ge 15n/16$. But $|N_1|\ge cn/2$, and so 
fewer than $n/16$ vertices in $Z_2$
have no neighbour in $N_1$, since $c>0$. Hence at least $7n/8$ 
vertices in $Z_2$ do have such a neighbour. 
Choose $B_1\subseteq N_1$ minimal such that
$B_1$ covers at least $5n/16$ vertices in $Z_2$.
Let $B_2$ be the set of vertices in $Z_2$ covered by $B_1$. Thus $5n/16\le |B_2|\le 3n/8$
from the minimality of $B$. Let $A_2=Z_2\setminus B_2$. Thus $A_2$ is anticomplete to $B_1$, 
and $|A_2|=|Z_2|-|B_2|\ge (15n/16-3n/8)=9n/16$. 

Let $A_1=N_1\setminus B_1$.
Since $|N_1|\ge cn/2$, the edge-density between $N_1,A_2$ is at least $c$. 
In particular there is a vertex $v_2\in A_1$
with at least $c|A_2|\ge 9cn/16\ge cn/2$ neighbours in $A_2$. 
(Note that $v_2\notin B_1$ since $B_1$ is 
anticomplete to $A_2$.) Let $N_2$ be the set of neighbours of $v_2$ in $A_2$.
Let $C_1$ be the set of vertices in $B_1$ adjacent to $v_2$, and let 
$D_2$ be the set of vertices in $B_2$ that have a neighbour in $B_1\setminus C_1$.
\\
\\
(1) {\em If $|D_2|\ge n/8$ then $G$ contains $C_5$.}
\\
\\
Assume that $|D_2|\ge n/8$. It follows that there is a set $D_2'\subseteq D_2$ of at least $n/16$
vertices that are nonadjacent to $v_2$.
The edge-density between $N_{2}$ and $D_2'$ is at least $c$, since $|N_{2}|\ge cn/2$,
and in particular some vertex $d_2\in D_2'$ has a neighbour $w\in N_{2}$. Since $d_2\in D_2'\subseteq D_2$, 
it is adjacent to
some vertex $d_1\in B_1$ that is nonadjacent to $v_2$; but then
$$d_1\d v_1\d v_2\d w\d d_2\d d_1$$
is an induced cycle of length $5$. (Note that $d_1$ is nonadjacent to $w$ since $B_1$ is anticomplete to $A_2$.)
This proves (1).

\bigskip

Let $Y_2=A_2\setminus N_2$; it follows that
$|Y_2|\ge |A_2|-n/16\ge n/2$. Since $|N_2|\ge cn/2$
the edge-density between $N_2, Y_2$ is at least $c$, and so some vertex $v_{3}\in N_2$ has at least $c|Y_2|\ge cn/2$
neighbours in $Y_2$. Let $N_3$ be the set of neighbours of $v_3$ in $Y_2$. Let $C_2$ be the set of vertices in $B_2$
with a neighbour in $C_1$.
\\
\\
(2) {\em If $|C_2|\ge 3n/16$ then $G$ contains $C_5$.}
\\
\\
Assume that $|C_2|\ge 3n/16$. It follows that there is a set $C_2'\subseteq C_2$ of at least $n/16$
vertices that are nonadjacent to both $v_2, v_3$. 
The edge-density between $N_3$ and $C_2'$ is at least $c$, since $|N_3|\ge cn/2$,
and in particular some vertex $c_2\in C_2'$ has a neighbour $w\in N_3$. Since $c_2\in C_2'\subseteq C_2$, it is adjacent to 
some vertex $c_1\in C_1$; but then 
$$c_1\d v_2\d v_3\d w\d c_2\d c_1$$
is an induced cycle of length $5$. (Note that $c_1$ is nonadjacent to $v_3, w$ since $B_1$
is anticomplete to $A_2$.)  This proves (2).

\bigskip

Since $B_1$ covers $B_2$, it follows that $C_2\cup D_2=B_2$, and since $|B_2|\ge 5n/16$,
the result follows from (1) and (2). This proves \ref{get5hole}.~\bbox


Next we apply a theorem of R\"odl~\cite{rodl}, the following.
($\overline{G}$ denotes the complement graph of $G$.)
\begin{thm}\label{rodl}
For every graph $H$ and all $d>0$ there exists $\delta>0$ such that for every $H$-free graph $G$,
there exists $X\subseteq V(G)$ with $|X|\ge \delta|V(G)|$ such that in one of $G[X]$, $\overline{G}[X]$,
every vertex in $X$ has degree at most $d|X|$.
\end{thm}

We deduce:

\begin{thm}\label{C5density}
There exists $\epsilon>0$ such that for all $c$ with $0\le c\le 1/2$, if $G$ is $C_5$-free with $n>1$
vertices, then there exist disjoint $A,B\subseteq V(G)$ with $|A|\ge \epsilon cn$ and $|B|\ge \epsilon n$,
such that the edge-density between $A,B$ is either less than $c$ or more than $1-c$.
\end{thm}
\Proof
Let $\delta$ satisfy \ref{rodl}, taking $d=1/20$ and $H=C_5$. Now let $\epsilon =\delta/16$,
and let $G$ be $C_5$-free with
$n>1$ vertices. 
Let $v$ be a vertex; then it has either at least $(n-1)/2$ neighbours or at least $(n-1)/2$ non-neighbours; and since
$(n-1)/2\ge \epsilon n$, we may assume that $1<\epsilon cn$, for otherwise the theorem holds taking $A=\{v\}$.
In particular $n>2\epsilon^{-1}\ge 32\delta^{-1}$.

By \ref{rodl}, there exists $X\subseteq V(G)$ with $|X|\ge \delta n$ such that 
every vertex of $J$ has degree at most $|V(J)|/20$, where $J$ is one of $G[X]$, $\overline{G}[X]$.
Since $|V(J)|\ge \delta n\ge 32$, it follows that
every vertex of $J$ has closed degree at most $|V(J)|/16$.
Since $C_5$ is isomorphic to its complement, $J$ is $C_5$-free, and so from \ref{get5hole},
there are two disjoint subsets $A,B\subseteq V(J)$ with $|A|\ge  c |V(J)|/2$ and $|B|\ge |V(J)|/16$,
such that the edge-density between $A,B$ in $J$ is at most $c$. Thus $|A|\ge c \delta n/2\ge \epsilon cn$ and 
$|B|\ge \delta n/16= \epsilon n$,
and the edge-density between $A,B$ in $G$ is either at most $c$ or at least $1-c$. This proves \ref{C5density}.~\bbox

It is possible to deduce versions of
\ref{EHpartial} from versions of R\"odl's theorem~\ref{rodl} directly, as follows. 
If we have $d,\delta$ satisfying
\ref{rodl}, then for any $n$, if we choose $k\le \min (\frac{1}{2d},\frac{\delta n}{2})$
then we can use Tur\'an's theorem to obtain a stable set or clique on $k$ vertices from the set of at
least $2k$ vertices with density at most $\frac{1}{2k}$ or at least $1-\frac{1}{2k}$ that \ref{rodl} gives us. 
This motivates trying to improve the bound in \ref{rodl}. 

\begin{itemize}
\item R\"odl's original proof of \ref{rodl} uses Szemer\'edi's regularity lemma and gives a
tower-type bound for $1/\delta$ in terms of $1/d$, which yields something worse than \ref{EHpartial}. 
\item In~\cite{FS}, a better bound of $\delta=2^{-15|V(H)|(\log (1/d)^2}$ in \ref{rodl} is proved, which implies
the bound of~\ref{EHpartial}.
\item It is conjectured that a polynomial dependence of $\delta$
on $d$ holds, and this would imply the Erd\H{o}s-Hajnal conjecture itself. 
\item For $H=C_5$ we can get mid-way between, and that provides a different route to proving \ref{mainthm}, as follows.
One can prove that for $H=C_5$ we may take 
$$\delta=2^{-O( \log(1/d)^2 / \log \log (1/d))}$$ 
in \ref{rodl} by 
appropriately adapting the
proof of \ref{rodl} in \cite{FS} using that we now know \ref{foxdensity} for $H=C_5$. This would imply \ref{mainthm}. 
But the details of the proof of this improved
bound for \ref{rodl} for $C_5$ are involved and similar to that of the proof of \ref{mainthm} given in the next section,
and we omit them for the sake of brevity.
\end{itemize}

\section{The proof of \ref{mainthm}.}

Now we use \ref{C5density} to prove \ref{mainthm}. Since the argument to come is rather heavy, and works 
just as well for any graph $H$ satisfying \ref{foxdensity} instead of $C_5$,  it might be wise to
present it in full generality. Thus, let us say a class of graphs $\mathcal{I}$ is {\em hereditary}
if every graph isomorphic to an induced subgraph of a member of the class also belongs to the class.
Let $\epsilon$ be as in \ref{C5density}, and let $\sigma>\log(\epsilon^{-1})$. Then for $c\le 1/2$, $c^\sigma\le \epsilon$, and so
by \ref{C5density}, if $G$ is $C_5$-free with $n\ge 2$ vertices, and $0\le c\le 1/2$, then there exist disjoint $A,B\subseteq V(G)$
with $|A|\ge c^\sigma n$ and $|B|\ge \epsilon n$, such that the edge-density between them is either at most $c$ or at least $1-c$.
Then \ref{mainthm} follows from \ref{C5density} and the following, applied to the hereditary class of all $C_5$-free graphs:

\begin{thm}\label{thm1}
Let $\mathcal{I}$ be a hereditary class of graphs, and let $\sigma\ge 0$ and $0\le \epsilon\le 1$ 
with the following property: for every graph $G\in \mathcal{I}$ with at least two vertices,
and all $c$ with $0\le c\le 1/2$, there are disjoint subsets $A,B\subseteq V(G)$ with $|A|\ge c^\sigma n$
and $|B|\ge \epsilon n$, such that the edge-density between $A,B$ is either at most $c$ or at least $1-c$, where $n=|V(G)|$.
Then there exists $\kappa>0$ such that 
$$\max(\alpha(G),\omega(G))\ge 2^{\kappa\sqrt{\log n \log \log n}}$$
for every $G\in \mathcal{I}$,
where $n=|V(G)|\ge 2$.
\end{thm}
\Proof Let us define $r(n)= \sqrt{\log n \log \log n}$ for $n\ge 2$, for typographical convenience.

A {\em cograph} is a graph
not containing a 4-vertex path.
Thus the disjoint union of two cographs is a cograph, and so is the complement of a cograph.
We prove \ref{thm1} by showing that $G$ contains a cograph with at least $2^{2\kappa r(n)}$ 
vertices. As cographs are perfect, there is a clique or independent set with $2^{\kappa r(n)}$ vertices (and so of the desired 
cardinality).

For a graph $G$, let $\phi(G)$ denote the maximum of $|V(H)|$ over all cographs $H$ contained in $G$.
For each real number $x\ge 0$, let $f(x)$ be the minimum
of $\phi(G)$, over all graphs $G\in \mathcal{I}$ with $|V(G)|= \lceil x\rceil$ (we may assume there is some such graph $G$,
or else the result is trivially true). It is easy to see that $f(x)$ is non-decreasing with $x$.

We may assume that $\sigma\ge 1$ (by increasing $\sigma$ if necessary).
Let $\mu = (32\sigma)^{-1/2}$.
Choose $n_0$ such that
$$\left \lfloor \frac{\sigma 2\mu r(n)-1}{\log(2/\epsilon)}\right \rfloor\ge \sqrt{\log n}$$
for all $n\ge n_0$, and also
such that
$\mu r(n_0)\ge 2$, and
$\log n_0\ge 4\sigma \mu r(n_0)$.
Choose $\kappa\le \mu/2$ such that $2\kappa r(n_0)\le 1$. We will show that $\kappa$ satisfies the theorem.
\\
\\
(1) {\em For all $n\ge 2$ and all $c$ with $0\le c\le 1/2$, either $f(n)\ge 1/(4c)$ or $f(n)\ge f(c^{\sigma}n/2)+f(\epsilon n/2)$.}
\\
\\
Let $G\in \mathcal{I}$
with $n\ge 2$ vertices, such that $\phi(G)=f(n)$. 
Since $G\in \mathcal{I}$, the hypothesis implies that there are 
disjoint sets $A,B\subseteq V(G)$ with $|A| \geq  c^\sigma n$ and $|B| \geq \epsilon n$ 
such that the edge-density between $A$ and $B$ is either at most $c$ or at least $1-c$. We suppose 
without loss of generality that this density is at most $c$ (in the other case, we apply the same argument 
to $\overline{G}$).

Let $A''$ be the set  of vertices in $A$ with at least $2c|B|$ neighbours in $B$. 
As the number of edges between $A,B$ is at least 
$2c|B||A''|$ and at most $c |A||B|$,
it follows that $|A''| \leq |A|/2$. 
Let $A'=A \setminus A''$; so $|A'|=|A|-|A''| \geq |A|/2$ and every vertex in $A'$ has at most $2c|B|$ 
neighbours in $B$. Since $G[A']\in \mathcal{I}$, it follows from the definition of $f$ that
$\phi(G[A'])\ge f(|A'|)$.
Let $A_0\subseteq A'$ induce a cograph, with $|A_0|=f(|A'|)$.

If $|A_0| > 1/(4c)$, then $f(n)=\phi(G)\ge |A_0|\ge 1/(4c)$ as required, so 
we may assume that $|A_0| \leq 1/(4c)$. 
Let $B'$ be those vertices in $B$ with no neighbours in $A_0$; so $|B'| \geq |B|-2c|B||A_0| \geq |B|/2$. 
Again from the definition of $f$, $\phi(G[B'])\ge f(|B'|)\ge f(\epsilon n/2)$. 
Since $A_0$ is anticomplete to $B'$, it follows that 
$$f(n)=\phi(G)\ge |A_0|+\phi(G[B'])\ge f(c^{\sigma}n/2)+f(\epsilon n/2).$$
This proves (1).
\\
\\
(2) {\em For all $n\ge 2$ and all $c$ with $0\le c\le 1/2$, if $\log n\ge \sigma\log(1/c)$ 
then either $f(n)\ge 1/(4c)$ or $f(n) \geq kf(c^{2\sigma}n)$, where
$$k=\left \lfloor \frac{\sigma \log(1/c)-1}{\log(2/\epsilon)}\right \rfloor.$$}

\smallskip

\noindent
We may assume that $f(n)< 1/(4c)$, and hence $f(n')<1/(4c)$ for all $n'\le n$.
From the definition of $k$, $k\log(2/\epsilon)\le \sigma \log(1/c)-1\le \log n -1$, and so
$n(\epsilon/2)^{k}\ge 2$. Hence we may 
recursively apply (1) $k$ times without violating the condition ``$n\ge 2$'' in (1); and we obtain
$$f(n) \geq f(c^{\sigma}n/2)+f(c^{\sigma}(\epsilon/2)n/2)+f(c^{\sigma}(\epsilon/2)^2n/2)+\cdots + f(c^{\sigma}(\epsilon/2)^{k}n/2).$$ 
Each of the $k+1$ terms on the right side is at least $f(c^{2\sigma}n)$, from the definition of $k$,
and so $f(n) \geq kf(c^{2\sigma}n)$. This proves (2).
\\
\\
(3) {\em For all $n\ge 2$ and all $c$ with $0\le c\le 1/2$, if $\log n\ge 2\sigma\log(1/c)$
and with $k$ as in {\rm (2)},
either $f(n)\ge 1/(4c)$ or $f(n) \geq k^j$, where
$$j = \left \lfloor \frac{\log n}{4\sigma\log (1/c)} \right \rfloor.$$}

\smallskip

\noindent Again, we may assume that $f(n)< 1/(4c)$, and hence $f(n')<1/(4c)$ for all $n'\le n$.
From the definition of $j$, $c^{2\sigma j}n\ge n^{1/2}$, and so $\log(c^{2\sigma j}n)\ge \frac12\log n\ge \sigma\log(1/c)$.
Moreover, $c^{2\sigma (j-1)}n\ge n^{1/2}c^{-2\sigma}\ge 2$ since $\sigma\ge 1$.
Hence we may apply (2) recursively $j$ times, and deduce that $f(n) \geq k^jf(c^{2\sigma j}n)\ge k^j$.
This proves (3).
\\
\\
(4) {\em Let $n\ge n_0$, and $c=2^{-2\mu r(n)}$. Then
\begin{itemize}
\item $c\le 1/2$;
\item $\log n\ge 4\sigma \mu r(n)$;
\item $k\ge \sqrt{\log n}$, where $k$ is as defined in {\rm (2)}; and 
\item $1/(4c)\ge 2^{\mu r(n)}$.
\end{itemize}
}
\smallskip
\noindent
We observe first that $c\le 1/2$ if $n\ge n_0$, since $\mu r(n_0)\ge 1$.
Also, $\log n_0\ge 4\sigma \mu r(n_0)$ from the choice of $n_0$,
and since $\frac{\log n}{r(n)}$ increases with $n$, it follows that
$\log n\ge 4\sigma \mu r(n)$ for $n\ge n_0$. But 
$4\sigma \mu r(n)= 2\sigma\log(1/c)$, and so 
$\log n\ge 2\sigma\log(1/c)$. This proves the second statement.
The third statement follows from the choice of $n_0$. 
For the final statement, we must check that
$\log(1/c)-2\ge \mu r(n_0)$, that is,
$2\mu r(n)\ge \mu r(n_0)+2$;
but $\mu r(n)\ge \mu r(n_0)$ since $n\ge n_0$,
and $\mu r(n)\ge 2$ from the definition of $n_0$.
This proves (4).
\\
\\
(5) {\em If $n\ge n_0$ then $f(n)\ge 2^{\mu r(n)}$.}
\\
\\
Let $c$ be as in (4) and let $n\ge n_0$. By 
the first two statements of (4); we may apply (3), and so either $f(n)\ge 1/(4c)$ or $f(n) \geq (\log n)^{j/2}$, by the third 
statement of (4).
In the first case, the claim follows from the final statement of (4), so we may assume that
$$f(n) \geq (\log n)^{j/2} \geq (\log n)^{(\log n)/(16\sigma\log(1/c))}  =  2^{(16\sigma \cdot 2\mu)^{-1}r(n)}.$$
As $\mu = (16\sigma \cdot 2\mu)^{-1}$
from the definition of $\mu$, this proves (5).

\bigskip

We recall that $\kappa\le \mu/2$ and $2\kappa r(n_0)\le 1$. We claim that $f(n)\ge 2^{2\kappa r(n)}$
for all $n\ge 2$. This is true if $n\le n_0$, because then $f(n)\ge 2\ge 2^{2\kappa r(n)}$; and
if $n>n_0$ then it follows from (5). This proves \ref{thm1}.~\bbox

\end{document}